\documentclass[11pt,a4paper]{amsart}
\usepackage[margin=1.1in]{geometry}
\usepackage{microtype}
\usepackage{hyphenat}
\usepackage{parskip}

\usepackage{amsmath,amssymb,amsthm,mathtools}
\usepackage{commath}
\usepackage{relsize}
\usepackage{microtype}
\usepackage{color}
\usepackage{graphicx}
\usepackage{array}
\usepackage{subfig}
\usepackage{wrapfig}
\usepackage{float}
\usepackage[normalem]{ulem}
\usepackage{outlines}
\usepackage{enumitem} % For clean list formatting
\usepackage{comment}
\usepackage{kantlipsum}

\usepackage{hyperref}
\usepackage{cleveref} 
\hypersetup{colorlinks=true, linkcolor=blue, citecolor=blue}

\crefname{enumi}{condition}{conditions}
\Crefname{enumi}{Condition}{Conditions}

\newcommand{\N}{\mathbb{N}}
\newcommand{\C}{\mathbb{C}}

\newcommand{\E}{\mathcal{E}}

\mathtoolsset{showonlyrefs=true}

\theoremstyle{plain}
\newtheorem{theorem}{Theorem}[section]
\newtheorem{corollary}[theorem]{Corollary}

\newtheorem{lemma}[theorem]{Lemma}

\theoremstyle{definition}
\newtheorem{remark}[theorem]{Remark}

\title[\bfseries Multiplicative functions additive on partitions of $2k$ nonzero squares]{Multiplicative functions additive on partitions of $2k$ nonzero squares}

\author{Jewel Mahajan}
\address{Department of Mathematical Sciences, Indian Institute of Science Education and Research (IISER) Berhampur, Ganjam, Odisha 760003, India}
\email{jewelmahajan421@gmail.com}

\subjclass[2020]{Primary 11A25; Secondary 11E25, 11N64}
\keywords{Multiplicative function, additive uniqueness, sums of two squares, Dubouis's theorem}

\begin{document}
\begin{abstract}
For a fixed integer $k \ge 3$, we study the multiplicative functions $f\colon\N\to\C$ satisfying
\[
   f\Bigl(\sum_{i=1}^{2k} x_i^2\Bigr) = \sum_{j=1}^{k} f\bigl(x_{2j-1}^2 + x_{2j}^2\bigr)
\]
for all positive integers $x_1,\dots,x_{2k}$. This extends a theorem of Park~\cite{ParkTwoSquares} on sums of two nonzero squares, which established the $k=2$ case. For $k=3$ and $k=4$, we prove that every such $f$ with $f(2)\neq 0$ is the identity function on $\N$. For $k \ge 5$, we show that such a function $f$ must be either the identity function on $\N$, or $f(n) = 0$ for all $n > 2k + 21$.
\end{abstract}

\maketitle

\section{Introduction}
An arithmetic function $f\colon \mathbb{N} \to \mathbb{C}$ is said to be \emph{multiplicative} if $f(1) = 1$ and $f(mn) = f(m)f(n)$ for all coprime $m, n \in \mathbb{N}$. Additionally, $f$ is \emph{completely multiplicative} if $f(mn) = f(m)f(n)$ holds for all $m, n \in \mathbb{N}$. Multiplicative functions are determined by their values at the prime powers. We call a function $f$ \textit{$k$-additive on a set $S\subseteq \N$}, if 
\[
    f(a_1 + \cdots + a_k) = f(a_1) + \cdots + f(a_k) \quad \text{for all } a_1, \dots, a_k \in S.
\]
A multiplicative function that is also $k$-additive for some $k \ge 2$ is forced to be the identity function in most cases. An early result in this direction is due to Spiro~\cite{Spiro1992}, who showed that if a multiplicative function $f$ is $2$-additive on the set of primes and $f(p_0)\neq 0$ for some prime $p_0$, then $f$ is the identity function (i.e., $f(n)=n$ for all $n \in \N$). This result was later extended to functions that are $k$-additive on the primes for $k=3$ by Fang~\cite{Fang2011}, and for all $k \ge 3$ by Dubickas and \v{S}arka~\cite{DubickasSarka2013}.
Following Spiro~\cite{Spiro1992}, we define a set $S \subseteq \mathbb{N}$ to be a \emph{$k$-additive uniqueness set} for a class $\mathcal{F}$ of arithmetic functions if $\mathcal{F}$ has exactly one element $f$ satisfying the $k$-additivity on $S$. The case $k=2$ corresponds to the standard definition of an \emph{additive uniqueness set}. The characterisation of additive uniqueness sets, along with the exceptional functions that arise in the absence of uniqueness, has become an area of active research.

Figurate numbers often serve as additive uniqueness sets for the class of multiplicative functions. Chung~\cite{Chung} demonstrated that uniqueness fails for the $2$-additive case on the set of squares by characterising all such multiplicative and completely multiplicative functions. In contrast, Park~\cite{ParkSquares} later showed that any multiplicative function that is $k$-additive on the set of squares for $k \ge 3$ must be the identity function. In other words, the set of squares is a $k$-additive uniqueness set for multiplicative functions whenever $k \ge 3$. Besides the set of squares, the sets of triangular and tetrahedral numbers were established as $2$-additive uniqueness sets for the class of multiplicative functions in the earlier work of Chung and Phong~\cite{ChungPhong}. Park~\cite{ParkTri} extended this result for $k\ge 3$, by showing that the set of triangular numbers is a $k$-additive uniqueness set for the collection of multiplicative functions. Kim et al.~\cite{Kim2018} showed that the set of all nonzero generalised pentagonal numbers is a $2$-additive uniqueness set for the collection of multiplicative
functions. This result was later extended to the $k$-additive case by Hasanalizade~\cite{Hasanalizade2022} for $k \ge 3$, where the proof is conditional on the Generalised Riemann Hypothesis for $k=3$ and unconditional for $k \ge 4$.

A related variant of these problems groups the summands into pairs, requiring additivity on sums of elements of $\E_2$ rather than on individual squares. Here
\[
   \E_2=\{\,a^2+b^2 : a,b\ge 1\,\}=\{2,5,8,10,13,17,18,20,25,26,\dots\}
\]
is the set of sums of two nonzero squares. Park~\cite{ParkTwoSquares} determined all multiplicative functions satisfying
\[
   f(a^2+b^2+c^2+d^2)=f(a^2+b^2)+f(c^2+d^2),\qquad a,b,c,d\ge 1.
\]
More precisely, Park showed that any such function $f$ is the identity function provided that $f(3)f(11) \neq 0$; otherwise, $f(n) = 0$ for all $n \ge 2$, except for $n \in \{3, 9, 11\}$. In other words, Park determined all functions that are $2$-additive on $\E_2$; the companion problem of classifying functions that are $2$-additive on $T_2$, the set of sums of two positive triangular numbers, was settled in~\cite{ParkTwoTri}. In each of these two-summand problems, the function is essentially forced to be the identity, apart from a small family of exceptional functions.

In this article, we deal with the higher analogue of~\cite{ParkTwoSquares} obtained by summing $k$ pairs. We study multiplicative functions $f$ satisfying
\begin{equation}\label{eq:main}
   f\Bigl(\sum_{i=1}^{2k}x_i^2\Bigr)=\sum_{j=1}^{k}f\bigl(x_{2j-1}^2+x_{2j}^2\bigr),
   \qquad x_1,\dots,x_{2k}\ge 1 .
\end{equation}
Because every element of $\E_2$ is of the form $x^2+y^2$ with $x,y\ge 1$,
\eqref{eq:main} is equivalent to
\begin{equation}\label{eq:reduced}
   f\Bigl(\sum_{i=1}^{k}A_i\Bigr)=\sum_{i=1}^{k}f(A_i)
   \qquad\text{for all } A_1,\dots,A_k\in\E_2 ,
\end{equation}
In other words, we classify the multiplicative functions that are $k$-additive on $\E_2$ for $k\ge 3$. The case $k=2$ is precisely \cite[Theorem~1]{ParkTwoSquares}, and we treat the remaining cases $k\ge 3$ in this article.

The change in analysis from $k=2$ to $k \ge 3$ comes from a structural  change in the representability of integers as sums of $2k$ nonzero squares. For $k=2$, the set of sums of $2k$ nonzero squares excludes an infinite set, such as $\{2\cdot 4^r, 6\cdot 4^r, 14\cdot 4^r \mid r\ge 0\}$, which complicates the analysis. By contrast, for $k = 3$ (resp., $k=4$), the set of sums of $2k$ nonzero squares contains all but finitely many integers, as every integer greater than or equal to $20$ (resp., $22$) is representable as a sum of six (resp., eight) nonzero squares, our approach relies, on which our approach relies. We utilise this property, which simplifies the classification of multiplicative functions beyond the case settled in~\cite{ParkTwoSquares}.

Throughout, the hypothesis $f(2)\neq 0$ serves as the nondegeneracy condition, similar to the role played by $f(p_0)\neq 0$ in Spiro's theorem. It cannot be omitted, because the function
$f$ with $f(1)=1$ and $f(n)=0$ for $n>1$ is multiplicative and satisfies
\eqref{eq:main} for every $k$. Our first result shows that
for $k=3$ and $k=4$, this is the only obstruction.
\begin{theorem}\label{thm:intro-k34}
Let $f$ be a multiplicative function that satisfies~\eqref{eq:main} with
$k\in\{3,4\}$. If $f(2)\neq 0$, then $f(n)=n$ for all positive integers $n$.
\end{theorem}
For $k \ge 5$, we establish that $f$ must eventually behave as either the identity or the zero function, regardless of the value of $f(2)$.
\begin{theorem}\label{thm:kge5-intro}
Let $k \ge 5$ and let $f$ be multiplicative and satisfy~\eqref{eq:main}. Then either $f(n) = n$ for all positive integers $n$, or $f(n) = 0$ for all $n > 2k + 21$.
\end{theorem}
In particular, if $f$ is a multiplicative function that satisfies~\eqref{eq:main} with $k\ge 5$, and $f(a)\neq 0$ for some $a>2k+21$, then $f$ is the identity function (see Corollary~\ref{cor:kge5-id}).

The proofs follow the same route as in~\cite{ParkTwoSquares}. For $k\in\{3,4\}$, one starts with the equivalent form~\eqref{eq:reduced} and uses small decompositions in $\E_2^{(k)}$, combined with multiplicativity and assumption $f(2)\neq0$, to evaluate $f$ on an initial segment of integers. The proof then continues via strong induction: since every integer beyond this segment is a sum of $2k$ nonzero squares, it can be expressed as a sum of $k$ elements of $\E_2$, all smaller than the integer itself. For $k\ge 5$, the base case is replaced by a simpler argument: by comparing the two decompositions $N=2+10+\sum C_i$ and $N+1=5+8+\sum C_i$ (where $C_i \in \E_2$), we show that $f(N+1)-f(N)$ is eventually constant for all $N>2k+21$, forcing $f$ to satisfy one of the forms in Theorem~\ref{thm:kge5-intro}. In all cases, we repeatedly use Corollary~\ref{cor:6sqand8sq}.

The remainder of the paper is organised as follows. Section~\ref{sec:prelim} collects the representability results we need, including Corollary~\ref{cor:6sqand8sq}.
Sections~\ref{sec:k3} and~\ref{sec:k4} prove Theorem~\ref{thm:intro-k34} by handling the cases $k=3$ and $k=4$ separately, and Section~\ref{sec:kge5} proves Theorem~\ref{thm:kge5-intro} and a subsequent corollary.

\section{Preliminaries}\label{sec:prelim}
Throughout, we write $x = f(2)$, $y = f(3)$, and $z = f(5)$. We will use the following result on the integers representable as a sum of a fixed number of nonzero squares.
\begin{lemma}[Dubouis~\cite{Dubouis1911}]\label{thm:dubouis}
Let $m \ge 4$. Every positive integer $n$ can be represented as a sum of $m$ nonzero squares, except for
\[
 \begin{cases} 
        \{1,3,5,9,11,17,29,41\}\cup\{\, 2\cdot 4^r, 6\cdot 4^r, 14\cdot 4^r \mid r \ge 0 \, \}, & \text{when } m = 4, \\
        \{1,2,3,4,6,7,9,10,12,15,18,33\}, & \text{when } m = 5, \\
        \{1,2,\dots,m-1\}\cup\{m+1,\,m+2,\,m+4,\,m+5,\,m+7,\,m+10,\,m+13\}, & \text{when } m \ge 6.
    \end{cases}
\]
In particular, for $m \ge 6$, every integer $n \ge m + 14$ is a sum of $m$ nonzero squares.
\end{lemma}
\begin{corollary}\label{cor:6sqand8sq}
Every integer $n \ge 20$ is a sum of $6$ nonzero squares, and every integer
$n \ge 22$ is a sum of $8$ nonzero squares.
\end{corollary}
Let $\E_{2}^{(k)} := \E_2 + \cdots + \E_2$ ($k$ summands). Since each element of $\E_2$ is a sum of two nonzero squares, the set of integers expressible as a sum of $2k$ nonzero squares is exactly $\E_{2}^{(k)}$. Hence, for $k \ge 3$, $\E_{2}^{(k)}$ contains all integers $n \ge 2k+14$, by Lemma~\ref{thm:dubouis}.

\begin{remark}
The consecutive integers $12 = 2+10$ and $13 = 5+8$ both lie in
$\E_2 + \E_2$, which is one of the key observations used in Section~\ref{sec:kge5}.
\end{remark}

\section{The case \texorpdfstring{$k = 3$}{k=3}}\label{sec:k3}
We begin by evaluating the function $f$ at the first few prime values, establishing a foundation for our later analysis of larger integers.
\begin{lemma}\label{lem:k3-base}
If $f$ is multiplicative, satisfies~\eqref{eq:reduced} with $k = 3$, and $f(2) \neq 0$, then $f(2) = 2$, $f(3) = 3$, $f(5) = 5$.
\end{lemma}

\begin{proof}
From the identity $3\cdot 2 = 2+2+2$, we have $y = 3$. From $ 2+2+8 = 2+5+5$, we get $2x + f(8) = x + 2z$, that is, $f(8) = 2z - x$. From $18 = 2+8+8$, we have $f(18) = x + 2f(8) = 4z - x$. On the other hand, $f(18) = x f(9)$ and $9 = 2+2+5$ gives $f(9) = 2x + z$. Therefore, we have
\begin{equation}\label{eq:k3-master}
    x(2x + z) = 4z - x.
\end{equation}
Thus the equality $2+10+10 = 2+2+18$ gives $x + 2xz = 2x + (4z - x)$, so $z(x - 2) = 0$.

If $z = 0$, then~\eqref{eq:k3-master} yields $2x^2 = -x$. Since $x \neq 0$, this implies $x = -\frac{1}{2}$, $f(8) = 2z-x=\frac{1}{2}$, and $f(10) = xz = 0$. From $20 = 5+5+10$, we have $f(20) = 0$, while $20 = 2+5+13$ yields $f(13) = \frac{1}{2}$. Therefore, $f(26) = x f(13) = -\frac{1}{4}$, while $26 = 8+8+10$ gives $f(26) = 1$, a contradiction. Hence $z\neq0$, which gives $x = 2$, and by~\eqref{eq:k3-master}, we have $z = 5$.
\end{proof}
\begin{theorem}\label{thm:k3}
Let $f$ be multiplicative and satisfy~\eqref{eq:main} with $k = 3$. If $f(2) \neq 0$, then $f(n) = n$ for all positive integers $n$.
\end{theorem}
\begin{proof}
For $k = 3$, since $f$ satisfies~\eqref{eq:main}, it also satisfies~\eqref{eq:reduced}. By Lemma~\ref{lem:k3-base}, we have $f(2)=2$, $f(3)=3$, $f(5)=5$, which gives $f(8)=2z-x=8$ and $f(9)=2x+z=9$. By multiplicativity, we have $f(6)=6$, $f(10)=10$, and $f(15)=15$.

Using the decompositions $12=2+5+5$, $14=2+2+10$, $17=2+5+10$, $18=2+8+8$, $20=5+5+10$, and $22=2+10+10$ in $\E_2^{(3)}$, we have
$f(n)=n$ for all $n \in \{12,14,17,18,20,22\}$.
Since $f(17)=17$, $17=2+2+13$ yields $f(13)=13$. By multiplicativity, $f(7)=f(14)/f(2)=7$, $f(11)=f(22)/f(2)=11$, and $f(4)=f(12)/f(3)=4$. 

Also, $32=2+10+20$ gives $f(32)=32$. The identities  $38=2+18+18$ and $48=8+8+32$ give $f(n)=n$ for $n \in \{38,48\}$. Therefore, $f(16)=f(48)/f(3)=16$ and $f(19)=f(38)/f(2)=19$. Combining everything with $f(1)=1$, we have $f(n)=n$ for all $1\leq n \leq 20$.

To conclude the proof, we now use the second principle of mathematical induction. Let $N \ge 20$ and $f(n)=n$ for all $n<N$. By Corollary~\ref{cor:6sqand8sq}, $N$ is a sum of $6$ nonzero squares, so $N = A_1 + A_2 + A_3$ with $A_i \in \E_2$. Since each $A_i\ge 2$, we have $A_i \leq N-4 < N$, and~\eqref{eq:reduced} gives $f(N) =\sum_{i=1}^{3}f(A_i) =\sum_{i=1}^{3} A_i =N$. 
\end{proof}

\section{The case \texorpdfstring{$k = 4$}{k=4}}\label{sec:k4}

\begin{lemma}\label{lem:k4-base}
If $f$ is multiplicative, satisfies~\eqref{eq:reduced} with $k = 4$, and $f(2) \neq 0$, then $f(n)=n$ for all $n\in \{2,3,4,5,6,7,8,10\}$.
\end{lemma}

\begin{proof}
By assumption, $x\neq 0$. From $8 = 2+2+2+2$, we have $f(8) = 4x$. From $14 = 2+2+2+8 = 2+2+5+5$, we have $7x = 2x + 2z$, so that $z = \frac{5}{2} x$, and $f(14) = 7x = x f(7)$ gives $f(7) = 7$. From $4\cdot 5 = 5+5+5+5$, we have $ f(4)z=4z $, which yields $f(4) = 4$ (since $z = \frac{5}{2} x\neq 0$).

From $8+8+8+8 = 2+10+10+10$, we have $ 16x=4f(8) = x + 3xz$, so that $15x=3xz$, yielding $z = 5$. From $z = 5$, together with  $z = \frac{5}{2} x$, we get $x=2$, hence $f(8) =4x= 8$. Finally, $3\cdot 8 = 2+2+10+10$, gives $f(24) = y f(8) = 8y= 2x + 2xz$. Since $x=2$ and $z=5$, we have $y = 3$. Furthermore, the multiplicativity of $f$ implies that $f(6) = 6$ and $f(10) = 10$.
\end{proof}

\begin{theorem}\label{thm:k4}
Let $f$ be multiplicative and satisfy~\eqref{eq:main} with $k = 4$. If $f(2) \neq 0$, then $f(n) = n$ for all positive integers $n$.
\end{theorem}
\begin{proof}
For $k = 4$, since $f$ satisfies~\eqref{eq:main}, it also satisfies~\eqref{eq:reduced}. By Lemma~\ref{lem:k4-base}, $f(n)=n$ for all $n \in \{2,3,4,5,6,7,8,10\}$.

Using the decompositions $11=2+2+2+5$, $14=2+2+2+8$, $16=2+2+2+10$, $17=2+2+5+8$, $19=2+2+5+10$, $20=5+5+5+5$, $22=2+2+8+10$, 
%$23=2+5+8+8$, $24=2+2+10+10$, $25=2+5+8+10$, 
and $26=2+8+8+8$ in $\E_2^{(4)}$, we have 
$f(n)=n$ for all $n \in \{11,14,16,17,19,20,22,26\}$.

By multiplicativity, $f(12)=12$, $f(15)=15$, $f(21)=21$, and $f(13)= f(26)/f(2) = 13$. The identity $4\cdot 9=8+8+10+10$ yields $f(9)=9$. Consequently, $f(18)=18$, by multiplicativity. Combining everything with $f(1)=1$, we have $f(n)=n$ for all $1\leq n \leq 22$.

Let $K \ge 22$ and $f(n)=n$ for all $n<K$. By Corollary~\ref{cor:6sqand8sq}, $K$ is a sum of $8$ nonzero squares, so $K = A_1 + \cdots + A_4$ with $A_i \in \E_2$. Since each $A_i\ge 2$, we have $A_i \leq K-6 < K$, and~\eqref{eq:reduced} gives $f(K) =\sum_{i=1}^{4}f(A_i) =\sum_{i=1}^{4} A_i =K$. The proof is now immediate from the second principle of mathematical induction.
\end{proof}
\begin{remark}
    The proof of Theorem~\ref{thm:intro-k34} follows immediately Theorems~\ref{thm:k3} and~\ref{thm:k4}.
\end{remark}
\section{The case \texorpdfstring{$k \ge 5$}{k>5}}\label{sec:kge5}
In this section, we provide a proof of Theorem~\ref{thm:kge5-intro}.
\begin{proof}[Proof of Theorem~\ref{thm:kge5-intro}]
Since $f$ satisfies~\eqref{eq:main}, $f$ also satisfies~\eqref{eq:reduced}. Let $N_0(k) \coloneqq 2k + 21$. Fix $N > N_0(k)$ and set $M = N - (2k+2)$. Since $N > 2k+21$, we have $M \ge 20$. Therefore, by Corollary~\ref{cor:6sqand8sq}, $M$ is a sum of $6$ non-zero squares, that is,
$M = C_3 + C_4 + C_5$ with $C_3, C_4, C_5 \in \E_2$. Setting $C_6 = \cdots = C_k = 2 \in \E_2$
(we assign nothing when $k = 5$), we obtain
\[
    \sum_{i=3}^{k} C_i = M + 2(k-5) = (N - 2k - 2) + (2k - 10) = N - 12,
    \qquad C_3,\dots,C_k \in \E_2 .
\]
Using $12 = 2 + 10$ and $13 = 5 + 8$, we have
\[
    N = 2 + 10 + \sum_{i=3}^k C_i, \qquad N + 1 = 5 + 8 + \sum_{i=3}^k C_i .
\]
Since $2,5,8,10 \in \E_{2}$, applying~\eqref{eq:reduced} to both yields
\[
    f(N) = f(2) + f(10) + \sum_{i=3}^k f(C_i), \qquad
    f(N + 1) = f(5) + f(8) + \sum_{i=3}^k f(C_i).
\]
Let $c = f(5) + f(8) - f(2) - f(10)$. Therefore, for all $N > N_0(k)$, we have
\[
    f(N+1) - f(N) = c.
\]
Hence, for all $n \ge N_{0}(k)+2$, the telescoping sum gives
\begin{align}
    &f(n)-f(N_0(k) + 1)\\
    =&\sum_{N=N_{0}(k)+1}^{n-1}(f(N+1) - f(N)) =\sum_{N=N_{0}(k)+1}^{n-1} c=c(n-N_{0}(k)-1).
\end{align}
Setting $d = f(N_0(k) + 1) - c(N_0(k) + 1)$, we have $f(n) = cn + d$ for all $n \ge  N_0(k)+2$.

For distinct primes $p, q$ with $p > q \ge N_0(k)+2$, multiplicativity of $f$ gives
\[
    cpq + d = (cp + d)(cq + d) = c^2 pq + cd(p + q) + d^2.
\]
Rearranging the above equation, we have
\[
    (c - c^2)pq - cd(p + q) + (d - d^2) = 0.
\]
Fix a prime $q_0\ge N_0(k)+2$. Then, for all primes $p>q_0$, we have
\[
    ((c - c^2)q_{0} - cd)p + (d - d^2) -cd q_{0} = 0.
\]
Since the above equality holds for infinitely many primes $p$, both coefficients must vanish, yielding
\[
    (c - c^2)q_{0} - cd=0, \qquad -cd q_{0} + (d - d^2)=0.
\]
Since $q_0\ge N_0(k)+2$ was arbitrary, we have
\[
    (c - c^2)q - cd=0, \qquad -cd q + (d - d^2)=0,
\]
for infinitely many primes $q\ge N_0(k)+2$. These further yields
\[
    c = c^2, \qquad cd = 0, \qquad d = d^2,
\]
so $(c, d) \in \{(0,0), (0,1), (1,0)\}$.

When $(c,d) = (0,1)$, we have $f(n) = 1$ for all $n \ge N_0(k)+2$. We now choose $A_1, \dots, A_k \in \E_2$ such that each exceeds $N_0(k)+2$. Then~\eqref{eq:reduced} gives \[1=f(\sum_{i=1}^{k}A_i)=\sum_{i=1}^{k} f(A_i) = k \ge 5,\] a contradiction. Hence $(c,d) \in \{(0,0),(1,0)\}$, that is, either $f(N) = N$ for all $N > N_0(k)$, or $f(N) = 0$ for all $N > N_0(k)$.

We now prove that if $f(N) = N$ for all $N > N_0(k)$, then $f$ is indeed the identity function on $\N$. Let $n\in \N$ be arbitrary. Pick a prime $p$ such that $p>\max (n,N_{0}(k))$. In particular, $pn>N_{0}(k)$ and $(p,n)=1$. By multiplicativity, $pn=f(pn)=f(p)f(n)=pf(n)$, which yields $f(n)=n$. Since $n\in \N$ was arbitrary, this proves $f$ is the identity function on $\N$.
\end{proof}

\begin{corollary}\label{cor:kge5-id}
Let $k \ge 5$ and let $f$ be multiplicative and satisfy~\eqref{eq:main}. If $f(a) \neq 0$ for some $a>2k+21$, then $f(n) = n$ for all positive integers $n$.
\end{corollary}
\begin{proof}
    The proof is immediate from Theorem~\ref{thm:kge5-intro}.
\end{proof}
\begin{remark}
The threshold $k = 5$ is sharp for this method and cannot be lowered, as we require the sumset of $2(k-2)$ nonzero squares to contain all but finitely many integers. By Theorem~\ref{thm:dubouis}, this holds for $2(k-2) \ge 5$ but fails at $k = 4$, where the sumset misses the infinite progressions $\{2 \cdot 4^r, 6 \cdot 4^r, 14 \cdot 4^r \mid r \ge 0\}$.
\end{remark}
\subsection*{Acknowledgements}
The author is supported by an institute postdoctoral fellowship from IISER Berhampur.
\bibliographystyle{plain}
\bibliography{references}
\end{document}